\documentclass[12pt]{article}

\title{The structure of typical clusters  \\
in large sparse random configurations} 

\author{Jean Bertoin\thanks{Laboratoire de Probabilit\'es, 
UPMC, 175 rue du Chevaleret, 75013 Paris; and DMA, ENS Paris, France. Email:
jean.bertoin@upmc.fr} \and 
Vladas Sidoravicius
\thanks{ IMPA,
Estr. Dona Castorina 110, 
Rio de Janeiro, Brasil; and CWI, Kruislaan 413, 1098 SJ
P.O. Box 94079, 1090 GB Amsterdam
The Netherlands. 
Email: vladas@impa.br} 
}\date{}

\usepackage{amsfonts,amsmath,amssymb,bbm}
\usepackage{setspace} 
\def\proof{\noindent{\bf Proof:}\hskip10pt}        
\def\QED{\hfill $\Box$}

\font\tenmath=msbm10 scaled 1200
\font\sevenmath=msbm7 scaled 1200
\font\Phiivemath=msbm5 scaled 1200
\newfam\mathfam \textfont\mathfam=\tenmath
\scriptfont\mathfam=\sevenmath \scriptscriptfont\mathfam=\Phiivemath

\def \\ { \cr }

\def \1{1 \mkern -6mu 1} 
\def\N{\mathbb{N}}
\def\E{\mathbb{E}}
\def\P{\mathbb{P}}
\def\GW{{\mathbb G \mathbb W}}

\def\D{{\bf D}}

\def \+{^{\rm (+)}}

\newtheorem{theorem}{Theorem}

\newtheorem{proposition}{Proposition}
\newtheorem{lemma}{Lemma}
\newtheorem{corollary}{Corollary}
\begin{document}

\maketitle

\begin{abstract} The initial purpose of this work is to provide a probabilistic explanation of a recent result on a version of Smoluchowski's coagulation equations in which the number of aggregations is limited.
The latter models the deterministic evolution of concentrations of particles in a medium where particles coalesce pairwise as time passes and each particle can only perform  a given number of aggregations. Under appropriate assumptions, the concentrations of particles converge as time tends to infinity to some measure which bears a striking resemblance with the distribution of the total population of a Galton-Watson process started from two ancestors. 

Roughly speaking, the configuration model is a stochastic construction which aims at producing a typical graph on a set of vertices with pre-described degrees.
Specifically, one attaches to each vertex a certain number of stubs, and then join pairwise the stubs uniformly at random to create edges between vertices.

In this work, we use the configuration model as the stochastic counterpart of
Smoluchowski's coagulation equations with limited aggregations. We establish a hydrodynamical type limit theorem for the empirical measure of the shapes of clusters in the configuration model when the number of vertices tends to $\infty$. The limit is given in terms of the distribution of a Galton-Watson process started with two ancestors.

  \end{abstract}

\begin{section}{Introduction}
The motivation for this work stems from a recent study of a deterministic model
for coagulation with limited number of aggregations. Specifically,  in  \cite{Bsolv},
one considers particles that are  determined by a pair of integers $(a,k)$ where
$k\geq 1$ represents the size and $a\geq 0$ the number of aggregations that the particle can perform. 
In the  model called symmetric, coagulations
 $$\{(a,k),(a',k')\}\,\longrightarrow \,(a+a'-2,k+k')$$ 
 occurs at rate 
$$aa'c_t(a,k)c_t(a',k')\,,$$ where $c_t(a,k)$ denotes the concentration of particles $(a,k)$ at time $t$ in the medium.  Analytically, this means that the evolution of concentrations is governed by the following variation of Smoluchowski's coagulation equations (cf. the survey by Aldous \cite{Ald2}):
\begin{eqnarray}\label{EQ1}
\frac{\rm d}{{\rm d}t} c_t(a,k)  =&\frac{1}{2}& \sum_{a'=1}^{a+1}\sum_{k'=1}^{k-1} a'(a-a'+2)
c_t(a',k')c_t(a-a'+2,k-k')\nonumber \\
&-&\sum_{a'=1}^{\infty}\sum_{k'=1}^{\infty}aa'c_t(a,k)c_t(a',k'),
\end{eqnarray}
where the first term in the right-hand side accounts for the creation of particles $(a,k)$
as the result of coagulations of pairs $\{(a',k'),(a-a'+2,k-k')\}$ and the second term
for the disappearance of particles $(a,k)$ after a coagulation with a particle $(a',k')$. 

One of the main results in \cite{Bsolv} is that under appropriate conditions on the initial data that we shall recall 
 later on,  the concentrations $c_t(a,k)$ have a limit as time $t$ tends to infinity which is given by
\begin{equation}\label{EQ2}
c_{\infty}(a,k)= {\bf 1}_{\{a=0\}}\frac{1}{k(k-1)}\nu^{*k}(k-2)\qquad \hbox{ for }
a\in\N \hbox{ and }k\geq 2\,.
\end{equation}
Here, $\nu$ is a certain probability measure on $\N$ with $\sum_{n=0}^{\infty} n\nu(n)\leq 1$ that depends on the initial data, and
$\nu^{*k}=\nu*\cdots * \nu$ denotes its $k$-th convolution power.  
The expression \eqref{EQ2} bears a striking resemblance with a special case of the celebrated formula due to Dwass \cite{Dwass} who established that the total population $T _2(\nu)$ generated by a (sub)-critical Galton-Watson branching process with reproduction law $\nu$ and started from two ancestors is given by
\begin{equation}\label{EQ3}
\P(T_2(\nu)=k)=\frac{2}{k}\nu^{*k}(k-2)\,,\qquad k\geq 2\,.
\end{equation}
This invites for a probabilistic explanation and provides the incentive  for the present work. 

Our approach for relating \eqref{EQ2} to \eqref{EQ3} stems from the fact that solutions to the classical  Smoluchowski's coagulation equations (without restriction on the number of aggregations) appear as the hydrodynamical limit of certain  stochastic coalescent models introduced by Marcus and Lushnikov. In some loose sense, the latter describe
the microscopic random dynamics of the particle system when the macroscopic evolution is governed by Smoluchowski's coagulation equations. This important feature has been established rigorously by Norris \cite{Norris}. We also refer to Section 5.2.1 in \cite{RFCP} for an elementary approach
in the special case of the multiplicative kernel, as the latter bears an obvious similarity  with \eqref{EQ1}.  On the other hand, it is well-known that the multiplicative coalescent  is naturally related to the size of the connected components in the random graph model of Erd\"os and R\'enyi, see in particular the remarkable paper by Aldous \cite{Ald1}.
This leads us to consider an extension of the random graph model where the sequence of degrees of vertices is given, and which is known  as the {\it configuration model}. 
Loosely speaking, the configuration model is constructed by an elementary stochastic algorithm which aims at producing a random graph on a set of vertices with pre-described degrees;
in general the resulting graph is not simple, in the sense that there may exist loops and
multiple edges. Typically,  a certain number of stubs is appended to each vertex, and one joins  pairwise the stubs uniformly at random to create edges between vertices.
This induces a natural partition 
of the set of vertices into clusters, i.e. connected components. 

Since its introduction independently by Bollob\'as \cite{Bollo} and Wormald \cite{Worm} (see also Bender and Canfield \cite{BC}), this model has been studied in the mathematical  literature by many authors. We refer e.g. to \cite{New} for an interesting review of applications of this and other random graph models to some real life network systems. The main known results chiefly concern asymptotics when the number of vertices is large and the empirical measure of the degrees of vertices converges. In particular, Molloy and Reed \cite{MR1} have determined the critical parameter for the existence of a giant component; see also \cite{MR2} and \cite{NSW}. In different directions,  van der Hofstad, Hooghiemstra and  co-authors \cite{EHH, HHM, HHZ} have made deep contributions to the study of distances between vertices in such random graphs, while Britton {\it et al.} \cite{BDML} used the configuration model to produce large random simple graphs with pre-described  asymptotic degree distribution.

If we neglect the appearance of multiple edges, loops or cycles 
which do not contribute to aggregation of clusters,  the configuration model may serve as a stochastic counterpart to the deterministic evolution of concentrations in the variant \eqref{EQ1}
of Smoluchowski's coagulation equations. This leads us to investigate the size of typical clusters, and more generally their combinatorial structures. Roughly speaking, the main result of this work is a hydrodynamical limit theorem for the empirical distribution of the shapes of clusters rooted at a generic stub. The limit is expressed in terms of a pair of Galton-Watson trees which are connected by an extra edge between the two roots. In particular, this yields the probabilistic explanation of the formal similarity
between the solution \eqref{EQ2}  and Dwass formula \eqref{EQ3}.

Let us now present some heuristics which are close to some of those that have already been used in the literature on configuration models  to relate the latter to Galton-Watson processes; see in particular \cite{MR1} and \cite{HHM}. Imagine that we pick a stub uniformly at random; the degree of the vertex to which this stub is appended has then the size-biased law of the degree of a typical vertex. We then pick a second stub uniformly at random to create the first edge. Informally, when the number of vertices is large, the degree of the vertex to which the second stub is appended has again the size-biased law and is essentially independent of the first. These two vertices should be viewed as the ancestors of two growing populations, where, by induction, individuals beget
independently and with a reproduction law given by the distribution of the {\it outer} degree
of a size-biased vertex. When the reproduction law is critical or sub-critical, the Galton-Watson process eventually becomes extinct, and extinction occurs before any loop, multiple edge or cycle arises in that cluster of the configuration model.  This suggests 
that the combinatorial structure of a typical (not too large) cluster could be described as a pair of independent Galton-Watson trees which are connected by an additional edge between the two roots. More precisely, the reproduction law should be given by the size-biased degree of a typical vertex, shifted by one unit, because the number of children corresponds to the outer-degree of the vertex.

The present work can be viewed as a companion to the recent paper \cite{BSV},  in which 
we also identify in terms of certain Galton-Watson trees the limiting empirical distribution of random structures that appear in a toy model for polymerization. More precisely,
we consider in \cite{BSV} a system of grabbing particles, where particles consist in monomers having a certain number of arms. Arms are activated successively uniformly at random, and each time an arm is activated, it grabs a particle uniformly at random amongst those which have not been previously grabbed and do not belong either to its own cluster. The main result of \cite{BSV} is that when  the  initial number of particles is large and the numbers of arms are given by i.i.d. random variables with mean less than $1$, then the empirical distribution of the shapes of polymers is closed to that induced by a Galton-Watson tree with a single ancestor and reproduction law given by the distribution of the number of arms of a typical monomer.

The plan of this work is as follows. The next section is devoted to preliminaries on configuration models, the combinatorial structure of planar rooted trees, and Galton-Watson processes. The emphasis is put on planar structures and their codings by the sequence of  degrees via breadth-first search.  The main result on the empirical distribution of the structures of rooted clusters in large random configurations is stated in Section 3 and then proved by explicit first and second moments estimates. Finally Section 4 is devoted to some applications. We shall point at certain invariance properties of Galton-Watson trees under random re-rooting,
and conclude by explaining the striking resemblance between the formulas \eqref{EQ2} and \eqref{EQ3}. 

\end{section}
\begin{section}{Preliminaries}

\subsection{Pairings, configurations and clusters}

The  aim of this work is to relate random configuration models to Galton-Watson trees, and as the latter have a natural planar structure, we shall introduce the former in planar setting which is tailored for our purposes. 
In this direction, we should imagine particles as planar star-shaped objects
consisting in a vertex to which a certain number of stubs are appended.

Formally, we consider some finite set  ${\mathcal V}$ of vertices and a map $d: {\mathcal V}\to \N^*$ where 
$d(v)$ represents the {\it degree} of the vertex $v$, that is number of stubs attached to $v$.  We denote by ${\mathcal S}={\mathcal S}({\mathcal V},d)$ the set of 
 stubs and  shall suppose for the sake of simplicity that the total number of stubs 
$$S:=\#{\mathcal S}=\sum_{v\in{\mathcal V}}d(v)$$ is even; otherwise we may always decide to add a new stub to some vertex (or to add a vertex with a single stub). 
We call a partition of ${\mathcal S}$ into $S/2$ pairs a {\it pairing of stubs}
and write $\Pi({\mathcal S})$ for the set of pairings of stubs. We first point at
the following elementary facts.

\begin{lemma}\label{L1}
{\rm (i)} The cardinal of $\Pi({\mathcal S})$ is given by
$$\#\Pi({\mathcal S})=\frac{S!}{(S/2)!} 2^{-S/2} = \prod_{i=1}^{S/2}(S-2i+1)\,.$$

\noindent {\rm (ii)} Consider a partition of ${\mathcal V}$ into two subsets 
${\mathcal V}_1,{\mathcal V}_2$  such that $S_1:=\sum_{v\in{\mathcal V}_1}d(v)$ and $S_2:=\sum_{v\in{\mathcal V}_2}d(v)$ are even numbers. Set ${\mathcal S}_1:={\mathcal S}({\mathcal V}_1,d)$
and ${\mathcal S}_2:={\mathcal S}({\mathcal V}_2,d)$. Then the map
$$(\pi_1,\pi_2)\longrightarrow \pi_1\sqcup \pi_2$$
is a bijection from $\Pi({\mathcal S}_1)\times  \Pi({\mathcal S}_2)$
to the subset of $\Pi({\mathcal S})$ consisting in pairings $\pi$ such that
there are no pairs $\{s_1,s_2\}$ in $\pi$ formed by a stub $s_1$ attached to a vertex in ${\mathcal V}_1$
and a stub $s_2$ attached to a vertex in ${\mathcal V}_1$.
\end{lemma}

\proof
Indeed, a generic pairing can be obtained by enumerating the stubs by $\{1,\ldots, S\}$  and then pairing the stubs according to the couples $(1,2), (3,4), \ldots,(S-1,S)$. There are $S!$ possible enumerations and the mapping is $(S/2)! 2^{S/2}$ on $1$, where $(S/2)!$ accounts for the number of permutations of the $S/2$ couples $(2i-1,2i)$, and $2^{S/2}$ for the number of
ways $S/2$ unordered pairs can be ordered into couples.
This establishes the first claim. The second is obvious. \QED

We then form {\it edges} $e=\{ v,v'\}$ with $v,v'\in{\mathcal V}$ by joining the tips of 
pairs of stubs $\{s,s'\}$, where $s$ (respectively, $s'$) is appended to $v$ (respectively, to $v'$).  We stress that an edge is unoriented, that it can be  a {\it loop} (i.e. the two vertices $v$ and $v'$ defining an edge may coincide), and 
that the same edge may appear by joining different pairs of stubs. 
Each pairing of stubs $\pi$ yields a {\it configuration} $\gamma(\pi)$, that is the family of the $S/2$ edges induced by the pairing. 
Note that there may be {\it multiple edges}; the same edge is repeated in $\gamma(\pi)$ as many times as it arises by joining different pairs of stubs in $\pi$.  We also stress that
the map $\pi\to \gamma(\pi)$ is not injective.

We view an edge which is not a loop as an elementary path connecting two different vertices, so a configuration $\gamma(\pi)$ on $({\mathcal V},d)$ naturally induces a partition of ${\mathcal V}$ into connected components. Endowing a given connected component with the restriction of  $\gamma(\pi)$ to the set of edges formed by pairs of vertices in that component, we obtain a  {\it cluster}.  
 
\subsection{Planar rooted trees and their structures}

Lemma \ref{L1}(ii) enables us to reduce the study of a given cluster to that of pairings
 $\pi\in \Pi({\mathcal S})$ such that the entire set of vertices ${\mathcal V}$ is connected for the configuration $\gamma(\pi)$. We shall therefore focus on that case in this section.
Recall that a {\it cycle} is a sequence of $\ell\geq 3$  distinct vertices, say $v_{1}, \ldots, v_{\ell}$, such that there exists an edge connecting $v_{j}$ and $v_{j+1}$ for every $j=1,\ldots, \ell-1$
 and also an edge connecting $v_{\ell}$ and $v_{1}$. 
 
A configuration $\gamma(\pi)$  that  connects ${\mathcal V}$ is called a {\it tree} if it contains no loops, no multiple edges, and no cycle. Note that this can occur only when
$S/2=\#{\mathcal V}-1$.
Because particles (i.e. vertices and the stubs that are appended) can be viewed as planar objects, we may think of tree-configurations as planar structures, 
in the sense that they can be represented in the plane in such a way that edges are line segments which do no cross, by attributing lengths to the edges in an appropriate manner. Throughout this section, we assume that $\#{\mathcal V}=k$ and that the configuration $\gamma(\pi)$ is a tree; in particular  $\gamma(\pi)$ consists in $k-1$ edges and $S=2(k-1)$.

To describe precisely the shape, that is the combinatorial structure, of a tree, we need to specify an origin and an orientation. For this, we distinguish a stub $s$ and call it the {\it root}.
This stub is  appended to a certain vertex $v$ that we use as the origin. Distinguishing ${s}$ also enables us to order all the stubs attached to $v$ by deciding that the first stub is ${s}$ and the next (if any) are ranked clockwise from that one.
Further, for every vertex $v'\neq v$ in that tree, we distinguish the stub appended to $v'$ that points at the origin $v$. 
This provides a natural order on the set of stubs  appended to any given vertex of the tree, and thus enables the use of {\it breadth-first search} to enumerate the vertices of the tree.

Specifically, set $s_1=s$ and $v_1=v$, define $s_{i+1}$ as the stub that is paired with the
$i$-th stub appended to  $v_1$ for $i=1,\ldots, d(v_1)$,
and write $v_{i+1}$ for vertex to which $s_{i+1}$ is appended.
We should think of $v_2, \ldots, v_{d(v_1)+1}$ as the children of $v_1$.
The stub $s_2$ is chosen as the first of the stubs appended to $v_2$, thus it is the unique stub pointing at the origin
and the other stubs attached to $v_2$ are ranked clockwise from $s_2$ and point at the children of $v_2$ (i.e. the vertices at distance $2$ from the origin $v_1$ and at distance $1$ from $v_2$). We denote these $d(v_2)-1$ children 
by $v_{d(v_1)+2}, \ldots, v_{d(v_1)+d(v_2)+1}$, and continue with the next
children $v_3, \ldots, v_{d(v_1)+1}$ of $v_1$ is an obvious way.
Then we proceed with indexing the third generation of vertices, in the order which
is naturally induced by the indexation of the second generation, and so on. See the figure below.

\begin{picture}(300,250)(-100,-30)
\put(10,10){\circle{20}}
\put(10,110){\circle{20}}
\put(60,60){\circle{20}}
\put(120,60){\circle{20}}
\put(170,60){\circle{20}}
\put(220,60){\circle{20}}
\put(120,10){\circle{20}}
\put(170,10){\circle{20}}
\put(120,110){\circle{20}}
\put(172,110){\circle{20}}
\put(70,170){\circle{20}}
\put(172,170){\circle{20}}
\put(222,170){\circle{20}}

\put(10,10){\makebox(0,0){$3$}}
\put(10,110){\makebox(0,0){$4$}}
\put(60,60){\makebox(0,0){$1$}}
\put(120,60){\makebox(0,0){$2 $}}
\put(170,60){\makebox(0,0){$6$}}
\put(220,60){\makebox(0,0){$10$}}
\put(120,10){\makebox(0,0){$7 $}}
\put(170,10){\makebox(0,0){$ 11$}}
\put(120,110){\makebox(0,0){$5$}}
\put(172,110){\makebox(0,0){$13 $}}
\put(70,170){\makebox(0,0){$8 $}}
\put(172,170){\makebox(0,0){$9$}}
\put(222,170){\makebox(0,0){$12$}}

 \put (18,18){\line(1,1){35}}
 \put (18,102){\line(1,-1){35}}
  \put (70,61){\line(1,0){10}}
   \put (70,59){\line(1,0){10}}
   \put(82,60){\makebox(0,0){$\triangleright$}}
  \put (84,60){\line(1,0){26}}
 \put (130,60){\line(1,0){30}}
\put (180,60){\line(1,0){30}}

 \put (120,50){\line(0,-1){30}}
\put (170,50){\line(0,-1){30}}
 \put (120,70){\line(0,1){30}}
 \put (172,160){\line(0,-1){39}}
\put (182,170){\line(1,0){30}}

 \put (126,119){\line(1,1){41}}
  \put (116,119){\line(-1,1){41}}

\end{picture}

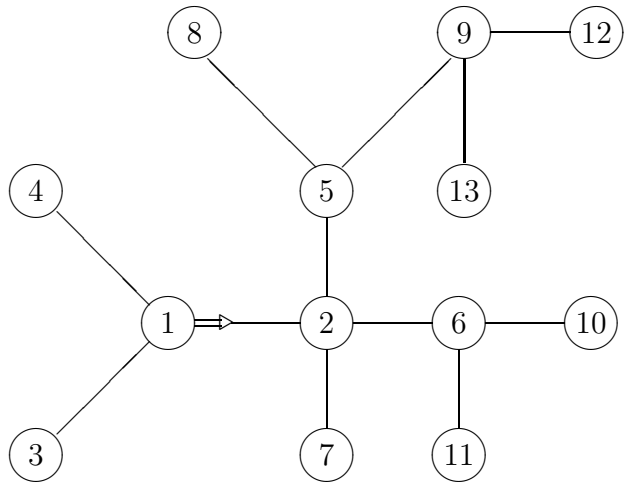
\centerline{\bf Figure 1 : \sl Enumeration by breadth-first search of the vertices of a planar tree
}
\centerline{\sl rooted at the stub $=\mkern -7mu \rhd$. The degree sequence is $(3,4,1,1,3,3,1,1,3,1,1,1,1)$.}

We write $d_i$ for the degree of the $i$-th vertex. We stress  that for $2\leq i\leq k$, 
the outer-degree of $v_i$, i.e. the number of stubs appended to $v_i$
that point away from the origin, is $d_i-1$. It is well-known that the sequence of degrees ${\bf d}=(d_1, \ldots, d_k)$ 
fulfills 
\begin{equation}\label{EQ4}
\min\{j\geq 1: d_1+\cdots+d_j=2(j-1)\}=k\,,
\end{equation}
and characterizes a unique planar rooted tree structure. 
Conversely, any finite sequence ${\bf d}=(d_1, \ldots, d_k)$  such that \eqref{EQ4} holds 
encodes a planar rooted tree structure with $k$ vertices.
We write $\D$ for the set of sequences ${\bf d}=(d_1, \ldots, d_k)$ which fulfill \eqref{EQ4}, where the lenght  $k\in\N^*$ is arbitrary, and  think
of the set $\D$ of sequences of degrees as the set of structures of planar rooted trees.
We refer for instance to Section 6.2 in Pitman \cite{PiSF} for details.

We now summarize this discussion, introducing first some terminology for convenience.  A bijection $\{1,\ldots, k\}\to {\mathcal V}$  
can be represented as a sequence ${\bf v}=(v_1,\ldots,v_k)$ of distinct vertices and will be referred to as an {\it enumeration} of ${\mathcal V}$. We also call
a map $\varsigma : {\mathcal V}\to {\mathcal S}$ 
that associates to each vertex $v\in{\mathcal V}$ a stub $s\in{\mathcal S}$ appended to that vertex a {\it selection of stubs}.
For every pairing $\pi\in\Pi({\mathcal S})$ such that the configuration $\gamma(\pi)$ on ${\mathcal V}$ is a tree and every choice of a distinguished stub $s\in{\mathcal S}$, breadth first search yields a 
unique enumeration  ${\bf v}=(v_1,\ldots,v_k)$ of ${\mathcal V}$
such that  the sequence $d({\bf v})=(d(v_1), \ldots, d(v_k))$ belongs to $\D$ (i.e. fulfills \eqref{EQ4}), and a unique  a selection of stubs $\varsigma$. The map
$$(\pi,s)\longrightarrow ({\bf v},\varsigma)$$
is bijective.  More precisely, we recover the pairing $\pi$ and the root stub $s$ by first constructing the planar rooted tree structure associated to ${\bf d}=(d(v_1), \ldots, d(v_k))$, and then placing the vertices
$v_1, \ldots, v_k$ on this structure in the order induced by the breadth first search.
The first stub appended to $v_1$ is $s=\varsigma(v)$,
and for every $i=2,\ldots,k$, $\varsigma(v_i)$ is the stub appended to $v_i$
which points at the origin  $v_1$. This determines the pairing $\pi$.

In order to record this analysis, it is convenient to introduce 
the multinomial coefficient
\begin{equation}\label{EQ5}
M({\mathcal V},d):= \left( 
\begin{matrix} k\\  \ell_1,\ldots,\ell_j\\
\end{matrix}\right)
=\frac{k!}{\ell_1!\cdots \ell_j!}\,,
\end{equation}
where $j$ is the number of different values, say $x_1, \ldots, x_j$, occurring in the family  $(d(v): v\in{\mathcal V})$, and $\ell_i$ the number of occurrences of the value $x_i$ in that family. For every structure ${\bf d}\in\D$, we say that ${\bf d}$ is {\it compatible} with $({\mathcal V},d)$ if there is at least an enumeration ${\bf v}=(v_1,\ldots,v_k)$ of ${\mathcal V}$ 
such that ${\bf d}=(d(v_1), \ldots, d(v_k))$, that is if and only if the sequence ${\bf d}$
takes the same values with the same multiplicity as the family $(d(v): v\in{\mathcal V})$.
The following statement should now be plain.

\begin{lemma} \label{L2}
 Suppose that $\#{\mathcal V}=k$ and $S=2(k-1)$. 
 Fix a rooted planar tree structure ${\bf d}=(d_1, \ldots, d_k)\in \D$.
If ${\bf d}$ is compatible with $({\mathcal V},d)$, then
 the number of pairs $(\pi,s)\in\Pi({\mathcal S})\times {\mathcal S}$
 for which the configuration $\gamma(\pi)$ is a tree with structure ${\bf d}$ when rooted at $s$ equals
 $$M({\mathcal V},d) \, \prod_{v\in{\mathcal V}}d(v).$$ 
 Otherwise (i.e. if ${\bf d}$ is not compatible), this number is $0$.
  \end{lemma}
  We stress that all the rooted planar tree structures which are compatible with
$({\mathcal V},d)$ are thus equally likely to occur 
if we choose the pair $(\pi,s)\in\Pi({\mathcal S})\times {\mathcal S}$ uniformly at random.
In the same vein, it may be also interesting to point at the following simple formula, 
even though it will not be used in this paper .

  \begin{proposition} \label{P1}
 Suppose that $\#{\mathcal V}=k$ and $S=2(k-1)$. 
 The number of pairings $\pi \in\Pi({\mathcal S})$
 for which the configuration $\gamma(\pi)$ is a tree, is
 $$(k-1)! \, \prod_{v\in{\mathcal V}}d(v),$$ 
\end{proposition}

 \proof 
 To establish the formula, we simply need to calculate the number of
enumerations ${\bf v}$ of ${\mathcal V}$ for which the sequence 
$(d(v_1), \ldots, d(v_k))$ corresponds to some rooted planar tree structure. 
Recall from the ballot theorem (see, e.g., Lemma 6.1 in \cite{PiSF}) that for each enumeration ${\bf v}$,  there is a unique
cyclic permutation $\sigma$ of $\{1, \ldots, k\}$ such that
$(d(v_{\sigma(1)}), \ldots, d(v_{\sigma(k)}))$ fulfills \eqref{EQ4}. 
This shows that this number is $(k-1)!$. 
\QED

\subsection {Galton-Watson trees with two ancestors}
We consider now a probability measure $\nu$ on $\N$ and associate to $\nu$  a  measure on $\D$ by
\begin{equation}\label{EQ6}
\GW^{\nu}_2({\bf d})=\prod_{i=1}^k\nu(d_i-1)\,,
\end{equation}
where ${\bf d}=(d_1, \ldots, d_k)$ denotes a generic rooted planar tree structure.

The measure $\GW^{\nu}_2$ has a simple interpretation in terms of Galton-Watson
branching processes, and is in fact a sub-probability.
More precisely, consider a Galton-Watson process with reproduction law $\nu$ and started from two ancestors. The process can be represented on the upper-half plane,
where the individuals at generation $\ell\in\N$ lie on horizontal line $y=\ell$,
in an order consistent with that of their respective parents, so that the edges 
(line-segments) linking parents to children do not cross each other. We further connect the 
two ancestors by an additional edge, and distinguish the stub attached to
the left-most ancestor that thus points at the right-most ancestor. 
This enables us to list individuals (vertices) by breadth first search just as in the preceding section. Observe that the degree of the left-most ancestor (i.e. the origin)
is distributed as $1+\xi$ where $\xi$ is a random variable with law $\nu$,
whereas the outer-degrees  of the other individuals (i.e. their numbers of children)
are given by independent copies of $\xi$.

The event when the total population is finite has probability one
if and only if  the reproduction law $\nu$ is critical or subcritical, 
i.e. $\sum_{i\in\N }i\nu(i)\leq 1$, and $\nu\neq \delta_1$. Restricting our attention to this event, the structure of this planar rooted tree is a random variable in $\D$ which has distribution $\GW^{\nu}_2$ and is defective in the supercritical case.

\vskip 2mm \noindent
 {\bf Remark.} In the case when $\nu$ is the Poisson distribution with parameter $
p\leq 1$, then it is easily checked that the law $\GW^{\nu}_2$ also describes the
law of the genealogical tree  of a Galton-Watson process with reproduction law $\nu$,
started from a {\it single} ancestor, and conditioned to have size at least 2.

 \end{section}

  \begin{section}{A limit theorem for typical rooted clusters}
 
 For each fixed integer $n$, we consider a set ${\mathcal V}_n$ of $n$ vertices and a function
 $d_n: {\mathcal V}_n\to \N^*$ that specifies the number of stubs appended to each vertex.
 We introduce the empirical distribution of the number of stubs
 $$\mu_n(i):=\frac{1}{n}\#\{v\in{\mathcal V}_n: d_n(v)=i\}\,,\qquad i\in\N^*\,.$$
 We write 
 $$S_n:=\sum_{v\in{\mathcal V}_n}d_n(v)=n\sum_{i=1}^{\infty} i\mu_n(i)$$
 for the total number of stubs, 
 assuming for simplicity that this quantity is even.
 Our basic assumption is that the limit
 \begin{equation} \label{EQ7}
 \lim_{n\to\infty} \mu_n(i):=\mu(i)
 \end{equation}
 exists for every $i\geq 1$, and that the average number of stubs
  $$n^{-1}S_n=\sum_{i=1}^{\infty} i\mu_n(i)$$
  converges as $n\to\infty$ to the first moment of $\mu$, i.e.
  \begin{equation} \label{EQ8}
 \lim_{n\to\infty} \sum_{i=1}^{\infty} i\mu_n(i)=\sum_{i=1}^{\infty} i\mu(i) :=m<\infty\,.
 \end{equation}
 
 We also denote by $\mu^*$ the probability measure on $\N^*$ which is obtained
 from $\mu$ by  size-biased sampling, that is
 $$\mu^*(i):= \frac{i\mu(i)}{m}\,,\qquad i\in\N^*\,.$$
  A standard application of Scheff\'e's lemma shows that \eqref{EQ7} can then be re-enforced to
 \begin{equation} \label{EQ9}
 \lim_{n\to\infty} i\mu_n(i)\frac{n}{S_n}=\mu^*(i)\qquad \hbox{ in $L^1(\N^*)$.}
 \end{equation}
 Finally, we introduce the probability measure $\nu$ on $\N$ 
 induced from $\mu^*$ by the shift $i\to i-1$ from  $\N^*$ to $ \N$, viz. 
 $$\nu(i)=\mu^*(i+1)\,, \qquad i\geq 0\,.$$

We write ${\mathcal S}_n$ for the set of stubs appended to vertices in ${\mathcal V}_n$. We  pick a pairing $\pi\in\Pi({\mathcal S}_n)$ uniformly at random, and denote by
$\Gamma_n:=\gamma(\pi)$ the resulting random configuration 
 on $({\mathcal V}_n, d_n)$. For every stub ${s}\in{\mathcal V}_n$, if the cluster  of $\Gamma_n$ which contains $s$ is a tree, then
 $T_{{s}}$ denotes the combinatorial structure which results from rooting that tree at  the stub ${s}$ (see Section 2.2), and otherwise, we decide that $T_{{s}}=\varnothing$. 
 
 We are interested in  the random variable
 $$\rho_n({\bf d}):=\frac{1}{S_n}\#\{{s}\in{\mathcal S}_n :   T_{{s}}= {\bf d}\}\,,\qquad {\bf d}\in\D$$
 which counts the proportion of stubs $s$ such that the cluster rooted at $s$ induced by $\Gamma_n$ is a tree with structure ${\bf d}$.  
 Similarly, we write
 $$\rho_n(\varnothing):=\frac{1}{S_n}\#\{{s}\in{\mathcal S}_n :   T_{{s}}= \varnothing\}$$
for the proportion of stubs $s$ such that the cluster containing $s$ induced by $\Gamma_n$ is not a tree.
 The collection $(\rho_n({\bf d}): {\bf d}\in\D)$
 should thus be viewed as a variant of the empirical measure of tree-clusters.  We 
 now able to state our main asymptotic result on large random configurations.

\begin{theorem} \label{T1} Assume that \eqref{EQ7} and \eqref{EQ8} hold. 
Then for every planar rooted tree configuration ${\bf d}\in\D$, the following limit holds in $L^2(\P)$ : 
$$\lim_{n\to\infty} \rho_n({\bf d})=
\GW^{\nu}_2({\bf d})\,.$$

If we further suppose that 
\begin{equation}\label{EQ10}
\sum_{i=1}^{\infty} i(i-2)\mu(i)\leq 0\,,
\end{equation}
and also exclude the degenerate case when $\mu$ is the Dirac point mass at $2$, then
$$\lim_{n\to\infty} \rho_n(\varnothing)=0\qquad \hbox{in $L^1(\P)$}.$$
 \end{theorem}

The condition \eqref{EQ10}  plays an important part for random configuration models.  According to a well-known result due to Molloy and Reed \cite{MR1}, when \eqref{EQ10} fails (assuming also some further technical conditions), then there is some constant $c>0$ such that with probability one, the random configuration $\Gamma_n$ contains almost surely a cluster of size at least $c n$ when $n$ is sufficiently large. The size of this giant component is estimated in \cite{MR2}. 
At the opposite, when \eqref{EQ10} holds
with a strict inequality (again assuming some further technical conditions), Molloy and Reed \cite{MR1} have shown that with probability one, the random configuration $\Gamma_n$ contains at most $n^{1/4}$ cycles and no cluster of size at least $n^{1/4}$ whenever $n$ is sufficiently large.  Note that  in the critical case when \eqref{EQ10} is an equality, 
Theorem \ref{T1} implies that the probability that there is  a cluster of size at least $\varepsilon n$ tends to $0$ for any $\varepsilon >0$, because 
$\GW^{\nu}_2$ is a probability measure on $\D$.
 
The proof of Theorem \ref{T1} relies on asymptotics for the first and second moments of $\rho_n({\bf d})$.
We first state:

\begin{lemma}\label{L3} We have
$$\lim_{n\to\infty}\E(\rho_n({\bf d})) =  \GW^{\nu}_2({\bf d})$$
for every ${\bf d}\in\D$.
\end{lemma}

\proof   Let  the structure  ${\bf d}=(d_1,\ldots, d_k)$ have size $k\geq 2$. 
 We write ${\mathcal V}'$ for a generic subset of ${\mathcal V}_n$ with $k$ vertices
 and ${\mathcal S}'$ for the set of stubs in ${\mathcal S}_n$ which are appended to vertices in ${\mathcal V}'$. 
 There are two cases. 
 
 If the unordered families of degrees $\{d(v'): v'\in{\mathcal V}'\}$
 and $\{d_i: 1\leq i \leq k\}$ do no coincide (recall that in such families, numbers are repeated according to their multiplicity), then there is no pairing of stubs for which the vertices of  ${\mathcal V}'$ are those of a  tree-cluster with structure ${\bf d}$ when properly rooted. We say that ${\mathcal V}'$ is {\it bad}.
 
 Otherwise, we say that ${\mathcal V}'$ is {\it good}. Introduce the set $G'$ of couples $(s,\pi)\in {\mathcal S}'\times \Pi({\mathcal S}_n)$ such that the cluster rooted at $s$ induced by the configuration $\gamma(\pi)$
is a tree whose set of vertices coincides with ${\mathcal V}'$ and has structure ${\bf d}$.
The cardinal of $G'$ can then be computed by combining Lemmas \ref{L1} and \ref{L2}. Since $\#{\mathcal S}'=2(k-1)$,
 one gets  
 \begin{equation}\label{EQ11} 
 \# G'=
 M({\bf d})\, \frac{(S_n-2(k-1))!}{(S_n/2-k+1)!} 2^{-S_n/2+k-1}\, \prod_{i=1}^k d_i\,,
 \end{equation}
where $M({\bf d})$ denotes the multinomial coefficient
$$
M({\bf d}):= \left( 
\begin{matrix} k\\  \ell_1,\ldots,\ell_j\\
\end{matrix}\right)
=\frac{k!}{\ell_1!\cdots \ell_j!}\,,
$$
with $j$  the number of different values in the sequence  ${\bf d}$ and $\ell_i$ the number of occurrences in ${\bf d}$ of the $i$-th value for $1\leq i \leq j$. 

 So it remains to estimate the number of good subsets ${\mathcal V}'$ with $k$ vertices, and for this we use a probabilistic argument. We sample   uniformly at random $k$ vertices in ${\mathcal V}_n$, say, $v_1,\ldots, v_k$, successively and without replacement.
 It should be plain from the hypothesis \eqref{EQ7} that when
 $n\to\infty$, the $k$-tuple of degrees $(d_n(v_1),\ldots d_n(v_k))$ converges in distribution to the $k$-tuple formed by i.i.d. variables with law $\mu$; in particular 
 the probability that $(d_n(v_1),\ldots, d_n(v_k))={\bf d}$ tends to
 $\prod_{i=1}^k \mu(d_i)$ as $n\to\infty$.
 We readily deduce that  the probability that the (unordered) family $\{d(v_1), \ldots, d(v_k)\}$ is good converges as $n\to\infty$ to
  $$\frac {k!}{M({\bf d})}\, \prod_{i=1}^k \mu(d_i)\,.$$
 As there are $n!/(n-k)!\sim n^k$  $k$-tuples of distinct vertices in ${\mathcal V}_n$ and as the map
 that transforms a $k$-tuple into an unordered set is $k!$ to $1$, we conclude that
 the number of good subsets in ${\mathcal V}_n$ is equivalent for large $n$ to
\begin{equation}\label{EQ12}
\frac {n^k}{M({\bf d})}\, \prod_{i=1}^k \mu(d_i)\,.
\end{equation}

  Recall from Lemma \ref{L1}(i) that 
  $$\#\left( {\mathcal S}_n\times \Pi({\mathcal S}_n)\right)  = S_n\frac{S_n!}{(S_n/2)!} 2^{-S_n/2}\,,$$
  and that  $d_i\mu(d_i)=m\nu(d_i-1)$, by definition.
Putting the pieces together, we find
 \begin{eqnarray*}
 \E(\rho_n({\bf d}))  &\sim& \frac {n^k}{ S_n }\, \frac{(S_n-2(k-1))! (S_n/2)! }{S_n! (S_n/2-k+1)!} 
 2^{k-1}\, \prod_{i=1}^k (d_i\mu(d_i))\\
  &\sim& \frac {n^k}{S_n}\, S_n^{-2(k-1)}\, (S_n/2)^{k-1} \, 2^{k-1}\, \prod_{i=1}^k (m\nu(d_i-1))\\
 &=& \frac {(nm)^k}{S_n^k}\,  \prod_{i=1}^k \nu(d_i-1)\,.
 \end{eqnarray*}
By \eqref{EQ6} and \eqref{EQ8}, this completes the proof. \QED

Lemma \ref{L3} essentially means that if we pick a stub $s$ uniformly at random in
${\mathcal S}_n$ and independently of the random configuration $\Gamma_n$, then 
the conditional distribution of the combinatorial structure of the random cluster rooted at $s$ given the event that this cluster is a tree, converges weakly as $n\to\infty$ to the Galton-Watson law $\GW^{\nu}_2$. Theorem \ref{T1} is a much stronger statement
that involves the empirical distribution of structures of clusters, and requires second moment estimates.
 
 \begin{lemma}\label{L4} We have
$$\lim_{n\to\infty}\E((\rho_n({\bf d}))^2) =  (\GW^{\nu}_2({\bf d}))^2$$
for every ${\bf d}\in\D$.\end{lemma}

\proof The argument is similar to that of Lemma \ref{L3}; in particular  we shall use the same notation and terminology. 
We start from the expression 
$$\E((\rho_n({\bf d}))^2)=S_n^{-2}\E\left(\#\{(s',s'')\in{\mathcal S}_n
\times {\mathcal S}_n :   T_{s'}=T_{s''}= {\bf d}\} \right)\,.$$

Let  ${\mathcal V}'$ and ${\mathcal V}''$ two generic subsets of ${\mathcal V}_n$, both with $k$ vertices,
 and write ${\mathcal S}'$ (respectively, ${\mathcal S}''$) for the set of stubs in ${\mathcal S}_n$ which are appended to vertices in ${\mathcal V}'$ (respectively, ${\mathcal V}''$). 
Note that for every stubs $s'\in{\mathcal S}'$ and $s''\in{\mathcal S}''$, the identity $T_{s'}=T_{s''}\neq \varnothing$
can occur only if ${\mathcal V}'$ and ${\mathcal V}''$ are both good and, either coincide or are disjoint.

We first consider the situation when ${\mathcal V}'={\mathcal V}''$.  Recall that for any  $s'\in{\mathcal S}'$, if $T_{s'}= {\bf d}$, then  ${\mathcal S}'={\mathcal S}''$has exactly $2(k-1)$ stubs.
 It follows from the proof of Lemma \ref{L3} that  the number of triplets $(s',s'',\pi)\in {\mathcal S}'\times {\mathcal S}''\times  \Pi({\mathcal S}_n)$ 
such that the cluster rooted at $s'$ induced by the configuration $\gamma(\pi)$
is a tree whose set of vertices coincides with ${\mathcal V}'$ and $T_{s'}=T_{s''}={\bf d}$,  
is bounded from above by 
 $$2(k-1)M({\bf d})\, \frac{(S_n-2(k-1))!}{(S_n/2-k+1)!} 2^{-S_n/2+k-1}\, \prod_{i=1}^k d_i\,;$$
 see \eqref{EQ11}. 
 Multiplying this by the number of good subsets ${\mathcal V}'$ in ${\mathcal V}_n$, that is approximatively by \eqref{EQ12}, we get a quantity which is small compared to 
 $$\#({\mathcal S}_n \times {\mathcal S}_n \times \Pi( {\mathcal S}_n))
 = S_n^2 \frac{S_n!}{(S_n/2)!} 2^{-S_n/2}
$$ 
when $n\to\infty$.
We conclude  that  in the evaluation of $\E((\rho_n({\bf d}))^2)$, the contribution of pairs of stubs $(s',s'')$ that belong to the same cluster becomes asymptotically negligible.

Next we consider the situation when ${\mathcal V}'$ and ${\mathcal V}''$ are good and disjoint. 
By calculations similar to those that yield \eqref{EQ12} in the proof of Lemma \ref{L3}, we get
that the number of good disjoint pairs of subsets $({\mathcal V}',{\mathcal V}'')$ in ${\mathcal V}_n$
is equivalent for large $n$ to
$$
\frac {n^{2k}}{M({\bf d})^2}\, \left(\prod_{i=1}^k \mu(d_i)\right)^2\,.
$$
It then follows from Lemmas \ref{L1} and \ref{L2}
that the number of triplets $(s',s'',\pi)\in {\mathcal S}'\times {\mathcal S}''\times  \Pi({\mathcal S}_n)$ 
such that  $T_{s'}=T_{s''}={\bf d}$ and the stubs $s'$ and $s''$ belong to disjoint clusters is close to 
$$ n^{2k}\, \frac{(S_n-4(k-1))! }{(S_n/2-2k+2)!} 
 2^{-S_n/2+2k-2}\, \left(\prod_{i=1}^k d_i\mu(d_i)\right)^2\,.$$
 
Putting the pieces together yields the estimate 
   \begin{eqnarray*}
 \E((\rho_n({\bf d}))^2  &\sim& \frac {n^{2k}}{S_n^2 }\, \frac{(S_n-4(k-1))! (S_n/2)! }{S_n! (S_n/2-2k+2)!} 
 2^{2k-2}\, \left( \prod_{i=1}^k d_i\mu(d_i)\right)^2\\
  &\sim& \frac {n^{2k}}{S_n^{2}}\, S_n^{-4(k-1)}\, (S_n/2)^{2k-2} \, 2^{2k-2}\, \left(\prod_{i=1}^k (m\nu(d_i-1))\right)^2\\
 &\sim&  \left(\prod_{i=1}^k \nu(d_i-1)\right)^2\,.
 \end{eqnarray*}
By \eqref{EQ6}, this shows our claim. \QED

We are now able to establish Theorem \ref{T1}.

\noindent{\bf Proof of Theorem \ref{T1}:}\hskip10pt
Combining  Lemmas \ref{L3} and \ref{L4}, we see that the variance of $\rho_n({\bf d})$ tends to $0$ as $n\to \infty$, which establishes the first claim. Assume now further that \eqref{EQ10} holds and that $\mu\neq \delta_2$. Equivalently, this means that the reproduction law $\nu$ of the Galton-Watson process is critical or sub-critical, and is not the Dirac mass at $1$. So extinction occurs a.s. and 
$$\sum_{{\bf d}\in\D}\GW^{\nu}_2({\bf d})=1\,.$$
As
$$\rho_n(\varnothing)=1-\sum_{{\bf d}\in\D}\rho_n({\bf d})\,,$$
Fatou lemma entails our second assertion. \QED

 \end{section}
\begin{section}{Some applications}
In this Section, we shall develop some consequences of our main result.
Recall that Theorem \ref{T1} implies that if one selects a stub uniformly at random
and independently of  a large random configuration that fulfills the conditions there,
then the structure of the cluster rooted at that stub has asymptotically the distribution 
$\GW^{\nu}_2$. This hints at an interesting property of invariance of such Galton-Watson trees under uniform random re-rooting. Recall the construction of the structure of a planar tree rooted at some stub as it has been presented in Section 2.2; Figure 2 below should explain better than words what is meant by re-rooting a
rooted planar tree at some stub.

\begin{corollary}\label{C1}  Suppose that $\nu$ is a critical or subcritical probability measure on $\N$ with $\nu\neq \delta_1$. Let $D$ be a random rooted planar tree structure
 with distribution $\GW^{\nu}_2$. Conditionally on $D$, select one of the
 $2(|D|-1)$ stubs of $D$ uniformly at random, and denote by $D'$ the new  structure obtained from $D$ by re-rooting at that stub. Then  $D'$ has again the law $\GW^{\nu}_2$. 
\end{corollary}

\begin{picture}(400,250)(-30,-30)

\put(80,10){\circle{20}}
\put(40,10){\circle{20}}
\put(40,50){\circle{20}}
\put(120,50){\circle{20}}
\put(80,90){\circle{20}}
\put(160,90){\circle{20}}
\put(40,130){\circle{20}}
\put(120,130){\circle{20}}
\put(120,170){\circle{20}}

\put (40,20){\line(0,1){20}}
\put (50,10){\line(1,0){20}}
 \put (85,18){\line(1,1){26}}
\put (115,58){\line(-1,1){11}}
\put (75,98){\line(-1,1){26}}
 \put (85,98){\line(1,1){26}}
  \put (125,58){\line(1,1){26}}
\put (120,140){\line(0,1){20}}

\put(80,10)
{\makebox(0,0){$2$}}
\put(40,10)
{\makebox(0,0){$1$}}
\put(120,50)
{\makebox(0,0){$4$}}
\put(40,50)
{\makebox(0,0){$3$}}
\put(80,90)
{\makebox(0,0){$5$}}
\put(160,90)
{\makebox(0,0){$6$}}
\put(40,130)
{\makebox(0,0){$7$}}
\put(120,130)
{\makebox(0,0){$8$}}
\put(120,170)
{\makebox(0,0){$9$}}

 \put (90,85){\line(1,-1){11}}
 \put (88,83){\line(1,-1){10}}
 \put(100,72){\makebox(0,0){$\vartriangleleft$}}


\put(320,10){\circle{20}}
\put(280,10){\circle{20}}
\put(280,50){\circle{20}}
\put(320,50){\circle{20}}
\put(360,50){\circle{20}}
\put(240,50){\circle{20}}
\put(280,90){\circle{20}}
\put(360,90){\circle{20}}
\put(360,130){\circle{20}}

 \put (290,9){\line(1,0){7}}
 \put (290,11){\line(1,0){7}}
 \put(300,10){\makebox(0,0){$\vartriangleright$}}
 \put (302,10){\line(1,0){8}}
 \put (320,20){\line(0,1){20}}
 \put (280,20){\line(0,1){20}}
  \put (280,60){\line(0,1){20}}
 \put (360,60){\line(0,1){20}}
  \put (360,100){\line(0,1){20}}
 \put (329,12){\line(1,1){29}}
 \put (270,13){\line(-1,1){28}}

\put(320,10){\makebox(0,0){$4$}}

\put(280,10){\makebox(0,0){$5$}}

\put(280,50){\makebox(0,0){$8$}}

\put(320,50){\makebox(0,0){$6$}}

\put(360,50){\makebox(0,0){$2$}}

\put(240,50){\makebox(0,0){$7$}}

\put(280,90){\makebox(0,0){$9$}}

\put(360,90){\makebox(0,0){$1$}}

\put(360,130){\makebox(0,0){$3$}}

\end{picture}

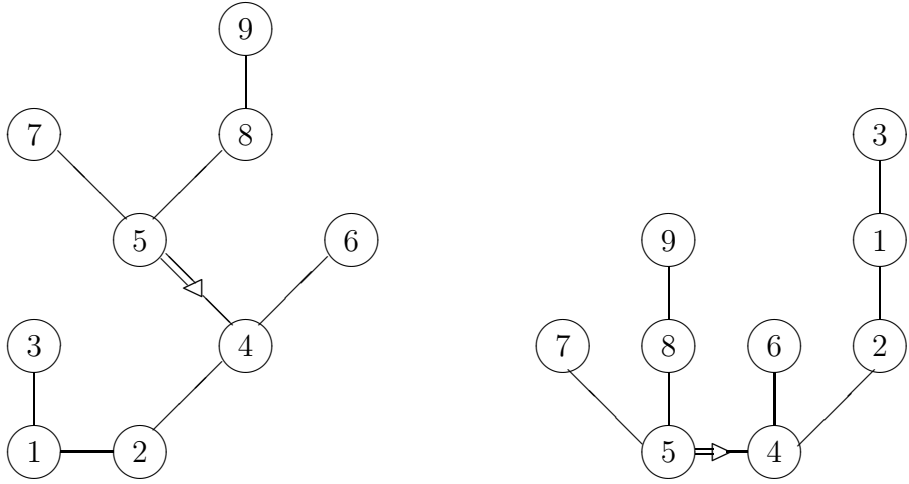
\centerline{\bf Figure 2 : \sl  Two genealogical  trees, both with two ancestors lying at the lowest level.}
\centerline{\sl  The left-most ancestor serves as the origin, the root-stub pointing at the right-most ancestor.}
\centerline{\sl The  tree on the right is the image of the tree on the left  by re-rooting at the stub $=\mkern -7mu \rhd$.
}
\centerline{\sl Vertices are labeled by breadth first order before re-rooting.}

\proof  Re-rooting has no effect on the degree of a vertex, so we only need to verify the statement for the conditional law of the Galton-Watson genealogical tree with two ancestors given the unordered  family of the degrees of vertices. 

Fix some unordered family, say $\Delta$, of $k$ positive integers (with possible repetitions), which add up to $2(k-1)$ and such that $\nu(\delta-1)>0$ for any integer $\delta$
in that family. Denote by $\D(\Delta)$ the subset of rooted planar tree structures 
 corresponding to some ordering of $\Delta$. 
We see from \eqref{EQ6} that the conditional law $\GW^{\nu}_2(\cdot \mid \D(\Delta))$
is simply the uniform distribution on $\D(\Delta)$. 

 Next consider the random configuration on a set $k$ vertices with degree family $\Delta$ that is induced by uniform random pairing, given that this configuration is a tree.
Then root the configuration using some stub that is picked independently and uniformly at random.  On the one hand, by construction, the law of the resulting combinatorial structure is obviously invariant by uniform random re-rooting. On the other hand, we see from Lemma \ref{L2} that it also coincides with the uniform distribution $\D(\Delta)$. This established our claim.  
\QED

We also refer to the recent work by Haas {\it et al.} \cite{HPW} and references therein
for a different property of invariance under uniform re-rooting for certain classes of random continuous trees.

It may be interesting to point also at the following avatar of Corollary \ref{C1}. 
A planar rooted tree is said {\it planted} if the degree of the origin, i.e. of the vertex to which the root-stub is appended, is $1$. In other words, the 
combinatorial structure ${\bf d}=(d_1,\ldots, d_k)$ fulfills $d_1=1$. So a planted
Galton-Watson tree describes the genealogy of a population where individuals beget independently with the same reproduction law, except the ancestor who has exactly one child. An easy consequence of Corollary \ref{C1} is that in the critical or sub-critical case, the structure of a planted Galton-Watson tree is statistically invariant under re-rooting at a leaf (i.e. a vertex with degree $1$) chosen uniformly at random.

We next turn our attention to some quantitative consequences of Theorem \ref{T1},
denoting for every $k\geq 2$ by $C_n(k)$ the number of clusters of size
$k$ in the random configuration $\Gamma_n$, i.e. the number of distinct connected components with $k$ vertices in the partition of ${\mathcal V}_n$ induced by $\Gamma_n$.

\begin{corollary}\label{C2} Assume that \eqref{EQ7}, \eqref{EQ8} and \eqref{EQ10} hold,
and exclude the case when $\mu=\delta_2$.
 We have
 $$\lim_{n\to\infty}\sum_{k=2}^{\infty}kÊ\E\left( \left|n^{-1}C_n(k)-\frac{m}{k(k-1)}\nu^{*k}(k-2)\right|\right)=0\,.$$
\end{corollary}

\proof Let us introduce first  for every $k\geq 2$ the subset  $\D_k$ of $\D$ consisting of structures of rooted planar trees ${\bf d}=(d_1,\ldots, d_k)$ of lenght $k$, and recall
that according to Dwass \cite{Dwass},
$$\GW^{\nu}_2(\D_k)=\frac{2}{k}\nu^{*k}(k-2)\,,$$
where $\nu^{*k}$ stands for the $k$-th convolution power of $\nu$.
As a tree of size $k$ has exactly $2(k-1)$ stubs and $\D_k$ is a finite set, we deduce from 
Theorem \ref{T1} that if we denote by $\tau_n(k)$ the number of clusters which are trees of size $k$, then
$$\lim_{n\to\infty} \frac{2(k-1)}{S_n}
\tau_n(k)=\frac{2}{k}\nu^{*k}(k-2)\,,$$
where the convergence takes place in $L^2(\P)$ and for every $k\geq 2$. 
Then we pick an arbitrary sequence of integers that tends to $\infty$, from which  we can excerpt  by a diagonal extraction procedure a subsequence such that with probability one,
$$\lim_{n\rightsquigarrow\infty} \frac{2(k-1)}{S_n}
\tau_n(k)=\frac{2}{k}\nu^{*k}(k-2) \qquad \hbox{for all $k\geq 2$,}$$
where the notation $n\rightsquigarrow \infty$ means that $n$ tends to infinity along that subsequence. 

Then observe that for each $n$, there is the obvious inequality
$$
\sum_{k\geq 2} 2(k-1)\tau_n(k)\leq S_n\,,$$
while 
$$\sum_{k\geq 2}\frac{2}{k}\nu^{*k}(k-2)=\sum_{k\geq 2}\GW^{\nu}_2(\D_k)=1\,,$$
since the reproduction law $\nu$ of the Galton-Watson process is critical or sub-critical and $\nu\neq \delta_1$.
A standard combination of Fatou and Scheff\'e lemmas entails that
$$\lim_{n\rightsquigarrow\infty} \sum_{k=2}^{\infty} \E\left(\left| \frac{2(k-1)}{S_n}
\tau_n(k)-\frac{2}{k}\nu^{*k}(k-2)\right| \right) =0\,.$$

Next, note that $\tau_n(k)\leq C_n(k)$ and $\sum_{k\geq 2} 2(k-1)C_n(k)\leq S_n$
as at least $2(k-1)$ distinct stubs are needed to connect  $k$ vertices. It follows that
\begin{eqnarray*}
& &\sum_{k=2}^{\infty} \E\left( \left|\frac{2(k-1)}{S_n}
(C_n(k)-\tau_n(k)) \right| \right)\\ &=&
\E\left(\sum_{k=2}^{\infty}  \frac{2(k-1)}{S_n}
C_n(k)\right)-\E\left(\sum_{k=2}^{\infty}  \frac{2(k-1)}{S_n}\tau_n(k)) \right)
\\
&\leq &
1-\E\left(\sum_{k=2}^{\infty}  \frac{2(k-1)}{S_n}\tau_n(k)) \right)\,,
\end{eqnarray*}
and we know from above that this quantity tends to $0$ as $n\rightsquigarrow\infty$.

This shows that
$$\lim_{n\rightsquigarrow\infty} \sum_{k=2}^{\infty} \E\left(\left| \frac{2(k-1)}{S_n}
C_n(k)-\frac{2}{k}\nu^{*k}(k-2)\right| \right) =0\,,$$
and since by the assumption \eqref{EQ8}, $S_n/n\to m$, we have thus proved that
$$\lim_{n\rightsquigarrow\infty} \sum_{k=2}^{\infty} k\E\left(\left| \frac{1}{n}
C_n(k)-\frac{m}{k(k-1)}\nu^{*k}(k-2)\right| \right) =0\,.$$
As the sequence of integers tending to infinity that we started from is arbitrary,
this establishes our claim. \QED

Corollary \ref{C2}  provides the explanation for the asymptotic behavior  \eqref{EQ2} 
 that motivated this work. Specifically, we know from
Theorem \ref{T1} that when the requirements \eqref{EQ7},  \eqref{EQ8}
and  \eqref{EQ10} are fulfilled, then, roughly speaking,  multiple edges, loops or cycles 
are rare. Roughly speaking, this means that almost all creation of edges correspond  to aggregations of clusters and thus enables us to view the configuration model as a stochastic microscopic version of  the terminal state of concentrations with a deterministic evolution governed by the variant \eqref{EQ1} of Smoluchowski's coagulation equations. In \cite{Bsolv}, one assumes that  initially all particles are monomers, i.e. consist in isolated vertices to which some stubs are appended. In the notation of the present work (beware that this differs from that in \cite{Bsolv}!), the initial concentration of particles with $i\geq 1$ stubs is 
$m^{-1}\mu(i)$, which is a finite measure on $\N^*$ with unit first moment.
In the framework of the random configuration model with $n$ vertices, this corresponds to assuming that particles live in a volume $mn$ and hence the
the initial concentration of monomers with $i$ stubs is given by
$$m^{-1}\mu_n(i)=
\frac{1}{mn}\#\{v\in{\mathcal V}_n: d_n(v)=i\}\,,\qquad k\in\N^*\,.$$
After the random pairing,  the  concentration of polymers with size $k$  (i.e. clusters with $k$ vertices)
is then $(mn)^{-1}C_n(k)$
and Corollary \ref{C2} shows that 
$$\lim_{n\to\infty} \frac{1}{mn}C_n(k):=c_{\infty}(0,k)= \frac{1}{k(k-1)}\nu^{*k}(k-2)\,.$$ One has thus
recovered  \eqref{EQ2}.

We also note that Corollary \ref{C2} solves a problem that has been addressed in Section II.C of \cite{NSW} by analytic and numerical technics. 
 \end{section}

\vskip 1cm
\noindent {\bf Acknowledgment.} We would like to thank Maria Eulalia Vares for 
stimulating discussions which have been at the origin of this work.

\end{document}